\documentclass{gtart_h}  

\def\ifplaintex{\expandafter\ifx\csname documentclass\endcsname\relax}

\ifplaintex 
\else
\fi

\expandafter\ifx\csname beginpicture\endcsname\relax
\expandafter\ifx\csname documentclass\endcsname\relax
\input pictex \else
\input prepictex \input pictex \input postpictex \fi\fi

\def\gt{{\mathsurround=0pt\it $\cal G\mskip-2mu$eometry \&\ 
$\cal T\!\!$opology}}        

\def\gtp{{\mathsurround=0pt\it $\cal G\mskip-2mu$eometry \&\ 
$\cal T\!\!$opology $\cal P\!$ublications}}  

\def\lognumber#1{\def\thelognumber{#1}}
\def\volumenumber#1{\def\thevolumenumber{#1}}
\def\papernumber#1{\def\thepapernumber{#1}}
\def\volumeyear#1{\def\thevolumeyear{#1}}

\def\pagenumbers#1#2{\def\startpage{#1}\def\finishpage{#2}}
\def\published#1{\def\publishdate{#1}}
\def\proposed#1{\def\theproposer{#1}}
\def\seconded#1{\def\theseconders{#1}}
\def\received#1{\def\receiveddate{#1}}
\def\revised#1{\def\reviseddate{#1}}
\def\accepted#1{\def\accepteddate{#1}}

\def\asciiaddress#1{\def\theasciiaddress{#1}}
\def\asciiemail#1{\def\theasciiemail{#1}}
\def\asciiurl#1{\def\theasciiurl{#1}}
\long\def\asciiabstract#1{\long\def\theasciiabstract{#1}}

\let\\\par\let\thelognumber\relax
\let\thevolumenumber\relax\let\thepapernumber\relax
\let\thevolumeyear\relax\let\thesamplenumber\relax\let\startpage\relax
\let\finishpage\relax\let\publishdate\relax\let\receiveddate\relax
\let\reviseddate\relax\let\accepteddate\relax\let\theasciititle\relax
\let\theasciiauthors\relax\let\theasciiaddress\relax
\let\theasciiabstract\relax
\let\theasciiemail\relax\let\theshortauthors\relax\let\theshorttitle\relax

\long\def\maketitlep{   

\count0=\startpage

\gt\hfill      
\beginpicture
\setcoordinatesystem units <0.33truein, 0.33truein> point at 2.2 0.9
\setplotsymbol ({$\cal G$})
\plotsymbolspacing=9truept
\circulararc 315 degrees from 0 1 center at 0 0
\setplotsymbol ({$\cal T$})
\circulararc 315 degrees from 1 -1 center at 1 0
\endpicture
\break
{\small\ifx\thesamplenumber\relax 
Volume \else Sample
\fi\thevolumenumber\ (\thevolumeyear)
\startpage--\finishpage\nl
Published: \publishdate}
\vglue 0.5truein plus 0.4fil minus 0.1truein

{\parskip=0pt\leftskip 0pt plus 1fil\def\\{\par\smallskip}{\ifplaintex\large
\else\Large\fi\bf\thetitle}\par\medskip}   

\vglue 0pt plus 0.1fil 

{\parskip=0pt\leftskip 0pt plus 1fil\def\\{\par}{\sc\theauthors}
\par\medskip}

\vglue 0pt plus 0.1fil 

{\small\parskip=0pt\let\newline\\
{\leftskip 0pt plus 1fil\def\\{\par}{\sl\theaddress}\par}
\expandafter\ifx\theemail\relax    
\relax\else\vglue 5pt plus 0.02fil minus 2pt\def\\{\stdspace{\rm 
and}\stdspace} 
\cl{Email:\stdspace\tt\theemail}\fi
\ifx\theurl\relax                  
\relax\else\vglue 5pt plus 0.02fil minus 2pt\def\\{\stdspace{\rm 
and}\stdspace}
\cl{URL:\stdspace\tt\theurl}\fi\par}

\vglue 7pt plus 0.3fil minus 3pt

{\bf Abstract}
\vglue 5pt plus 0.1fil minus 2pt

\theabstract

\vglue 7pt plus 0.3fil minus 3pt

{\bf AMS Classification numbers}\quad Primary:\quad \theprimaryclass

Secondary:\quad \thesecondaryclass

\vglue 5pt plus 0.3fil minus 2pt

{\bf Keywords:}\quad \thekeywords

\vglue 10pt plus 0.5fil minus 5pt

{\small  Proposed: \theproposer\hfill Received: \receiveddate\nl
Seconded: \theseconders\hfill 
\ifx\reviseddate\relax                         
Accepted: \accepteddate                        
\else
Revised: \reviseddate                          
\fi}
\eject
}       

\let\maketitlepage\maketitlep
\let\maketitle\maketitlepage

\font\phead=cmsl9 scaled 950
\font=cmsl9 scaled 1050
\font\pnum=cmbx10 scaled 913
\font\lnum=cmbx10 
\font\pfoot=cmsl9 scaled 950
\font=cmsl9 scaled 1050
\ifplaintex
\headline{\vbox to 0pt{\vskip -4.5mm\line{\small\phead\ifnum
\count0=\startpage ISSN 1364-0380 (on line)
1465-3060 (printed) \hfill {\pnum\folio}\else\ifodd\count0\def\\{ }%
\ifx\theshorttitle\relax\thetitle\else\theshorttitle\fi\hfill{\pnum\folio}
\else\def\\{ and }{\pnum\folio}\hfill\ifx\theshortauthors\relax\theauthors
\else\theshortauthors\fi\fi\fi}\vss}}
\footline{\vbox to 0pt{\vglue 0mm\line{\small\pfoot\ifnum\count0=\startpage
\copyright\ \gtp\hfill\else
\gt, Volume \thevolumenumber\ (\thevolumeyear)\hfill\fi}\vss
}}
\else
\makeatletter
\def\@oddhead{{\small\lhead\ifnum\count0=\startpage ISSN 1364-0380 (on line)
1465-3060 (printed) \hfill {\lnum\number\count0}\else\ifodd\count0
\def\\{ }\ifx\theshorttitle\relax \thetitle \else\theshorttitle\fi\hfill
{\lnum\number\count0}\else\def\\{ and }{\lnum\number\count0}
\hfill\ifx\theshortauthors\relax 
\theauthors\else\theshortauthors\fi\fi\fi}}\def\@evenhead{\@oddhead}
\def\@oddfoot{\small\lfoot\ifnum\count0=\startpage\copyright\ \gtp\hfill\else
\gt, Volume \thevolumenumber\ (\thevolumeyear)\hfill\fi}
\def\@evenfoot{\@oddfoot}
\makeatother
\fi

\newwrite\gtoutfile
\long\gdef\makeheadfile{  
{\def\\{, }\def\s{ }
\immediate\openout\gtoutfile head.xxx
\immediate\write\gtoutfile{Proxy-for: \ifx\theasciiauthors\relax
\theauthors\else\theasciiauthors\fi\s<\ifx\theasciiemail\relax\theemail\else\theasciiemail\fi>}
\immediate\write\gtoutfile{\noexpand\\}
\immediate\write\gtoutfile{Authors: \ifx\theasciiauthors\relax
\theauthors\else\theasciiauthors\fi}
{\def\\{ }\immediate\write\gtoutfile{Title: \ifx\theasciititle\relax
\thetitle\else\theasciititle\fi}}
\immediate\write\gtoutfile{Subj-class: GT or SG or MG etc}
\immediate\write\gtoutfile{MSC-class: \theprimaryclass\ifx\thesecondaryclass\relax\else, \thesecondaryclass\fi}
\immediate\write\gtoutfile{Journal-ref: Geom. Topol. \thevolumenumber
(\thevolumeyear) \startpage-\finishpage}
\immediate\write\gtoutfile{Comments: Published by Geometry and Topology at}
\immediate\write\gtoutfile{\s\s http://www.maths.warwick.ac.uk/gt/GTVol\thevolumenumber/paper\thepapernumber.abs.html}
\immediate\write\gtoutfile{\noexpand\\}
\immediate\write\gtoutfile{}
\ifx\theasciiabstract\relax
\immediate\write\gtoutfile{\theabstract}\else
\immediate\write\gtoutfile{\theasciiabstract}\fi
\immediate\write\gtoutfile{}
\immediate\write\gtoutfile{\noexpand\\}
\immediate\write\gtoutfile{}
\immediate\closeout\gtoutfile}}  

\def\maketitlepage{\maketitlep\makeheadfile}
\let\maketitle\maketitlepage

\lognumber{330}
\received{5 June 2003}
\volumenumber{9}\papernumber{2}\volumeyear{2005}
\pagenumbers{95}{119}   
\revised{20 December 2004}
\published{22 December 2004}
\accepted{29 September 2004}
\proposed{Martin Bridson}
\seconded{Cameron Gordon, Joan Birman}

\usepackage{amsmath,amssymb,epsfig,psfrag}

\newcommand{\HH}{{\mathbb{H}}}
\newcommand{\NN}{{\mathbb{N}}}
\newcommand{\RR}{{\mathbb{R}}}

\newcommand{\ZZ}{{\mathbb{Z}}}

\renewcommand{\setminus}{{\smallsetminus}}

\newcommand{\proj}{{\operatorname{proj}}}

\newcommand{\st}{\mid}
\newcommand{\from}{\co} 

\newcommand{\bdy}{{\partial}} 

\newcommand{\ML}{{\mathcal{ML}}} 
\newcommand{\PML}{{\mathcal{PML}}}

\theoremstyle{plain}
\newtheorem{theorem}{Theorem}[section]
\newtheorem{corollary}[theorem]{Corollary}
\newtheorem{lemma}[theorem]{Lemma}

\theoremstyle{definition}
\newtheorem{define}[theorem]{Definition}
\newtheorem{remark}[theorem]{Remark}
\newtheorem{claim}[theorem]{Claim}

\newcommand{\Kobayashi}{{\ensuremath{\mathbf{K}}}}
\newcommand{\Masur}{{\ensuremath{\mathbf{M}}}}
\newcommand{\Hempel}{{\ensuremath{\mathbf{H}}}}
\newcommand{\Bounded}{{\ensuremath{\mathbf{B}}}}
\newcommand{\Lin}{{\ensuremath{\mathbf{Lin}}}}
\newcommand{\Unb}{{\ensuremath{\mathbf{Unb}}}}

\newcommand{\MinLam}{{\ensuremath{\operatorname{MinLam}}}} 
\newcommand{\avedisplace}{{\alpha}}
\newcommand{\displace}{{\alpha}}

\begin{document}

\title{Distances of Heegaard splittings}

\author{Aaron Abrams\\Saul Schleimer}
\address{Department of Mathematics, Emory University\\Atlanta, 
Georgia 30322, USA\\{\rm and}\\Department of Mathematics, Rutgers 
University\\ Piscataway, New Jersey 08854, USA\\\medskip\tt
\\{\rm Email:\qua}\mailto{abrams@mathcs.emory.edu}{\rm\qua 
and\qua}\mailto{saulsch@math.rutgers.edu}\\\medskip
\\{\rm URL's:\qua}\url{http://www.mathcs.emory.edu/~abrams}\\
\url{http://www.math.rutgers.edu/~saulsch}}

\asciiaddress{Department of Mathematics, Emory University\\Atlanta, 
Georgia 30322, USA\\and\\Department of Mathematics, Rutgers 
University\\ Piscataway, New Jersey 08854, USA}

\asciiemail{abrams@mathcs.emory.edu, saulsch@math.rutgers.edu}

\asciiurl{http://www.mathcs.emory.edu/ abrams,
http://www.math.rutgers.edu/ saulsch}

\begin{abstract}
J~Hempel [Topology, 2001] showed that the set of {\em distances} of
the Heegaard splittings $(S, \mathcal{V}, h^n(\mathcal{V}))$ is
unbounded, as long as the stable and unstable laminations of $h$ avoid
the closure of $\mathcal{V} \subset \PML(S)$.  Here $h$ is a
pseudo-Anosov homeomorphism of a surface $S$ while $\mathcal{V}$ is
the set of isotopy classes of simple closed curves in $S$ bounding
essential disks in a fixed handlebody.

With the same hypothesis we show the distance of the splitting $(S,
\mathcal{V}, h^n(\mathcal{V}))$ grows linearly with $n$, answering a
question of A~Casson.  In addition we prove the converse of Hempel's
theorem.  Our method is to study the action of $h$ on the curve
complex associated to $S$.  We rely heavily on the result, due to
H~Masur and Y~Minsky [Invent.\ Math.\ 1999], that the curve complex
is Gromov hyperbolic.
\end{abstract}

\asciiabstract{%
J Hempel [Topology, 2001] showed that the set of distances of the
Heegaard splittings (S,V, h^n(V)) is unbounded, as long as the stable
and unstable laminations of h avoid the closure of V in PML(S).  Here
h is a pseudo-Anosov homeomorphism of a surface S while V is the set
of isotopy classes of simple closed curves in S bounding essential
disks in a fixed handlebody.

With the same hypothesis we show the distance of the splitting (S,V,
h^n(V)) grows linearly with n, answering a question of A Casson.
In addition we prove the converse of Hempel's theorem.  Our method is
to study the action of h on the curve complex associated to S.  We
rely heavily on the result, due to H Masur and Y Minsky
[Invent. Math. 1999], that the curve complex is Gromov
hyperbolic.}


\maketitlepage


\section{Introduction}

J~Hempel~\cite{Hempel01} introduced a new measure of the
complexity of a Heegaard splitting called the {\em distance} of the
splitting.  This definition is a conscious extension of A~Casson and
C~Gordon's notion of {\em strong
irreducibility}~\cite{CassonGordon87}.  Indeed, a Heegaard splitting
is
\begin{itemize}
\item
reducible if and only if its distance is 0,
\item
weakly reducible if and only if its distance is at most 1, and
\item
strongly irreducible if and only if its distance is at least 2.
\end{itemize}
Hempel's distance derives its name from the curve complex; the
distance of a splitting is exactly the distance between the two
handlebodies, thought of as subsets of the curve complex associated to
the splitting surface.  In particular the distance does not rely on
any particular diagram for the splitting.

Casson and Gordon produce examples of strongly irreducible splittings
by taking an existing splitting of $S^3$ and altering the gluing map
by high powers of a certain Dehn twist.  It is clear that the Dehn
twist must be carefully chosen; a Dehn twist about a curve which
bounds a disk in one of the two handlebodies leaves the splitting
unchanged.

Hempel obtains examples of high distance splittings using a
construction due to T~Kobayashi~\cite{Kobayashi88b}.  Instead of
$S^3$ he begins with the double of a handlebody $V$.  Instead of a
Dehn twist he iterates a certain pseudo-Anosov map $h$ on $S = \bdy
V$.  Analyzing the dynamics of $h$ acting on the space $\PML(S)$ of
projective measured laminations reveals that the set of distances,
obtained by altering the original gluing by $h^n$, is unbounded.  We
sketch Hempel's proof in Section~\ref{Sec:HempelsArgument}.  Again,
one must be careful when choosing $h$; if $h$ extends over the
handlebody then the splitting remains unchanged.  

We remark that the second author has proved that each fixed 3--manifold
has a bound on the distances of its Heegaard splittings.  Thus the
splittings provided by Hempel's theorem must represent infinitely many
different 3--manifolds.

This paper is part of an ongoing program to understand handlebodies
and Heegaard splittings from the point of view of the curve complex.
The fundamental ingredient underlying our approach is the result of
H~Masur and Y~Minsky~\cite{MasurMinsky99} that the curve complex is
$\delta$--hyperbolic.  This allows us to study the dynamics of $h$
acting on the curve complex.  We are thus able to both strengthen
Hempel's theorem and to prove a converse.

\begin{theorem}
\label{Thm:Main}
Fix a handlebody $V$ with genus at least two and set $S = \bdy V$.
Fix also a pseudo-Anosov map $h \from S \to S$.  Let $\mathcal{V}$ be
the set of isotopy classes of simple closed curves in $S$ which bound
disks in $V$.  Let $\avedisplace(h)$ denote the average
displacement of $h$.  The following are equivalent:

\begin{list}
{}{\setlength{\rightmargin}{\leftmargin}}
\item[\Kobayashi:] The stable and unstable laminations of $h$ are each
  of full type with respect to some pants decompositions of $V$.

\item[\Masur:] In $\PML(S)$, the stable and unstable laminations of
  $h$ are contained in the Masur domain of $V$.

\item[\Hempel:] In $\PML(S)$, the stable and unstable laminations of
  $h$ are not contained in the closure of ${\mathcal V}$.

\item[\Bounded:] In the curve complex, the projection of $\mathcal{V}$
  onto an invariant axis for $h$ has finite diameter.  

\item[\Lin:] There is a constant $K > 0$ so that for any $n \in \NN$
  the distance of the Heegaard splitting $(S, \mathcal{V},
  h^n(\mathcal{V}))$ lies between $n \cdot \avedisplace(h) - K$ and $n
  \cdot \avedisplace(h) + K$.  

\item[\Unb:] The set of distances of $\{(S, \mathcal{V},
  h^n(\mathcal{V})) \st n \in \NN\}$ is unbounded.
\end{list} 
\end{theorem}


The terms {\em average displacement}, {\em full type}, {\em Masur
domain}, and {\em invariant axis} are defined in
Definitions~\ref{Def:AveDisp},~\ref{Def:FullType} (or the
paper~\cite{Kobayashi88b}),~\ref{Def:MasurDomain} (or the
paper~\cite{Masur86}), and~\ref{Def:Axis}, respectively.

Hempel~\cite{Hempel01} proved that \Hempel\ implies \Unb.
Section~\ref{Sec:HempelsArgument} gives Hempel's definition of
distance as well as a sketch of his proof.

In this paper we introduce the condition \Bounded\ and prove in
Section~\ref{Sec:HIffB} its equivalence with \Hempel.  This
strengthens Hempel's theorem as \Bounded\ (and hence \Hempel) implies
\Lin.  (Note that \Lin\ implies \Unb\ by the fact, also contained
in~\cite{MasurMinsky99}, that $\avedisplace(h) > 0$.)  To prove this
we first develop several tools from $\delta$--hyperbolic geometry in
Sections~\ref{Sec:MetricSpaces} through~\ref{Sec:BoundedProj}.  These
arguments are written out carefully to emphasize their synthetic
nature, in particular the fact that they apply to spaces (such as the
curve complex) which are not locally compact.  In
Section~\ref{Sec:BImpliesLin}, we apply these tools to the curve
complex, as allowed by Masur and Minsky's theorem, to obtain
Corollary~\ref{Cor:BImpliesLin}, \Bounded\ implies \Lin.  The argument
in Section~\ref{Sec:HIffB} shows that \Bounded\ is an accurate
translation of \Hempel\ to the geometric language of the curve
complex.  This relies on E~Klarreich's
characterization~\cite{Klarreich99} of the Gromov boundary of the
complex of curves as the space of unmeasured, minimal laminations.

In order to prove the converse, namely that \Unb\ implies \Bounded, we
use the more recent theorem of Masur and Minsky~\cite{MasurMinsky03}
that handlebody sets are quasi-convex subsets of the curve complex.
This is carried out in Section~\ref{Sec:UnbImpliesB}.  The equivalence
of \Hempel, \Masur, and \Kobayashi\ is established in the final two
sections.
\medskip

{\bf Acknowledgements}\qua We thank Howard Masur for many interesting
conversations and for showing us the proof of
Lemma~\ref{Lem:HImpliesK'}.  AA was supported in part by NSF grant
DMS-0089927;  SS was supported in part by NSF grant DMS-0102069.

\section{Hempel's argument}
\label{Sec:HempelsArgument}

Beginning with a few definitions, this section states Hempel's
theorem and sketches a proof.  

\subsection{Terminology for Heegaard splittings}
A {\em handlebody} is a compact three-manifold which is homeomorphic
to a closed regular neighborhood of a finite, polygonal, connected
graph embedded in $\RR^3$.  The {\em genus} of the handlebody is the
genus of its boundary.  A properly embedded disk $D$ in a handlebody
$V$ is {\em essential} if $\bdy D$ is not null-homotopic in $\bdy V$.

Fix handlebodies $V$ and $W$ of the same genus.  Let $S = \bdy V$.
Glue $V$ and $W$ together via a homeomorphism $f \from \bdy V \to \bdy
W$.  We will consistently use $\mathcal{V}$ to denote the {\em
handlebody set}: the set of (isotopy classes of) curves in $S$
which bound essential disks in $V$.  Let $\mathcal{W}$ denote the set
of curves in $S$ which, after gluing, bound essential disks in $W$.
Then the data $(S, \mathcal{V}, \mathcal{W})$ give a {\em Heegaard
splitting}.  Note that a Heegaard splitting specifies a closed
orientable three-manifold.  The surface $S$ is referred to, in other
literature, as the {\em Heegaard splitting surface}.

\begin{define}
\label{Def:Distance}
The {\em distance} of the splitting $(S, \mathcal{V}, \mathcal{W})$ 
(see~\cite{Hempel01})
is the smallest $n \in \NN$ such that there are
$n + 1$ essential simple
closed curves $\alpha_i \subset S$ with the following properties:
\begin{itemize}
\item
$\alpha_0 \in \mathcal{V}$ and $\alpha_n \in \mathcal{W}$ and
\item 
$\alpha_i \cap \alpha_{i + 1} = \emptyset$ for $i = 0, \ldots, n-1$. 
\end{itemize}
\end{define}

Now suppose that $h$ is a homeomorphism of $S = \bdy V$.  Then the
set of curves $h(\mathcal{V})$ also defines a handlebody and $(S,
\mathcal{V}, h(\mathcal{V}))$ also specifies a Heegaard splitting.

\subsection{Hempel's theorem}
As above fix a handlebody $V$ of genus at least two 
and set $S = \bdy V$.  We will freely use
known facts about $\PML(S)$, the projectivization of the space of
measured laminations, and about the mapping class group of
$S$.  (See~\cite{FLP91} or~\cite{Kapovich01} for extensive discussion
of these objects.)  

For convenience of notation we occasionally blur the distinction
between an object and its isotopy class.  That said, let
$\mathcal{C}^0(S)$ be the set of isotopy classes of essential simple
closed curves in $S$.  Let $\mathcal{V} \subset \mathcal{C}^0(S)$ be
the set of curves which bound essential disks in the handlebody
$V$.  Fix also a pseudo-Anosov map $h \from S \to S$.

\begin{define}
\label{Def:H}
Given $V$, $S$, and $h$ as above, we say that \Hempel\ holds if the
stable and unstable laminations of $h$, $\mathcal{L}^{\pm}(h)$, are
not contained in the closure of $\mathcal{V}$ (considered as a subset
of $\PML(S)$).
\end{define}

\begin{define}
\label{Def:Unb}
Given $V$, $S$, and $h$ as above, we say that \Unb\ holds if the set
of distances of $\{ (S, \mathcal{V}, h^n(\mathcal{V})) \st n \in \NN \}$
is unbounded.
\end{define}

The following theorem of Hempel's~\cite{Hempel01} provided the first
proof that high distance splittings exist.  

\begin{theorem}[Hempel]
\label{Thm:Hempel}
Suppose a handlebody $V$ with $S = \bdy V$ and a pseudo-Anosov map $h
\from S \to S$ are given.  Then \Hempel\ implies \Unb.
\end{theorem}

\begin{remark}
Hempel cites Kobayashi~\cite{Kobayashi88b} as the framer of the
proof sketched below.  However, Kobayashi used a slightly different 
hypothesis; see Section~\ref{Sec:KIffH}.
\end{remark}

\begin{proof}[Sketch of proof of Theorem~\ref{Thm:Hempel}]
Suppose the distance of $(S, \mathcal{V}, h^j(\mathcal{V}))$ is
bounded by some $M \in \NN$, for all $j \in \NN$.  Then, for every $j$,
there is a set of essential curves $\{\alpha_i^j\}_{i = 0}^M \subset
\mathcal{C}^0(S)$ such that
\begin{itemize}
\item
$\alpha_0^j \in \mathcal{V}$ and $\alpha_M^j \in h^j(\mathcal{V})$ and
\item 
$\alpha_i^j \cap \alpha_{i + 1}^j = \emptyset$ for $i = 0, 1, \ldots,
M-1$.
\end{itemize}
Thus there is a curve $\beta^j \in \mathcal{V}$ such that
$h^j(\beta^j) = \alpha_M^j$.  By \Hempel\ the unstable
lamination $\mathcal{L}^-(h)$ is not contained in the closure of
$\mathcal{V}$ (taken in $\PML(S)$).  It follows that the $\beta^j$'s
avoid an open neighborhood of $\mathcal{L}^-(h)$.  Thus, the curves
$h^j(\beta^j) = \alpha_M^j$ converge to $\mathcal{L}^+(h)$ as a
sequence in $\PML(S)$.

Recall that $\PML(S)$ is compact. 
Inductively pass to subsequences exactly $M$ times to ensure that the
$i^{\rm th}$ sequence $\{\alpha_i^j\}_{j \in \NN}$ also converges
in $\PML(S)$, for $i = M - 1, M - 2, \ldots, 0$.  Denote the limit of
$\{\alpha_i^j\}_{j \in \NN}$ by $\mathcal{L}_i$.  Thus
$\mathcal{L}_M = \mathcal{L}^+(h)$ while $\mathcal{L}_0$ lies in the
closure of $\mathcal{V}$.  In particular $\mathcal{L}_M \ne
\mathcal{L}_0$ by \Hempel.

Recall that the geometric intersection number, $\iota(\cdot, \cdot)$,
extends to a continuous function $\ML(S) \times \ML(S) \to \RR_+$.
As $\alpha_i^j \cap \alpha_{i + 1}^j = \emptyset$ we have
$\iota(\mathcal{L}_i, \mathcal{L}_{i + 1}) = 0$, abusing notation
slightly.  As $\mathcal{L}^+(h)$ is minimal (the lamination contains
no closed leaf and all complementary regions are disks) 
and uniquely ergodic (all transverse measures are projectively
equivalent), $\mathcal{L}_M$ and $\mathcal{L}_{M - 1}$ must be the
same point of $\PML(S)$.
Inductively, $\mathcal{L}_i = \mathcal{L}_{i-1}$ in $\PML(S)$ which
implies that $\mathcal{L}_M = \mathcal{L}_0$.  This is a
contradiction.
\end{proof}

\subsection{Distance grows linearly}
The primary goal of this paper, then, is to show that Hempel's
hypothesis \Hempel\ implies a seemingly stronger assertion: the
distance grows linearly with the number of iterates of $h$, up to a
bounded additive constant.  This is the condition \Lin.  As indicated
in the introduction we do this by studying the action of $h$ on the
{\em complex of curves}, $\mathcal{C}(S)$.

\section{Metric spaces}
\label{Sec:MetricSpaces}

This section briefly states the facts we need about
$\delta$--hyperbolic spaces and their isometries.  For a deeper
discussion consult Gromov~\cite{Gromov87} or Bridson and
Haefliger~\cite{Bridson99}.

\subsection{Basic terminology}

Let $(X, d_X)$ be a metric space.  If $Y$ and $Z$ are subsets of $X$
define $d_X(Y, Z) = \inf\{d_X(y, z) \st y \in Y, z \in Z \}$.  

An {\em arc} in $X$ is a continuous map $L \from [a, b] \to X$ where
$[a, b]$ is a closed connected subset of the real numbers, $\RR$.  The
arc $L$ is {\em geodesic} if $|b' - a'| = d_X(L(a'), L(b'))$ for every
finite subinterval $[a', b'] \subset [a, b]$.  An arc $L \from [a,b]
\to X$ {\em connects} two points $u, v \in X$ if $[a, b]$ is a finite
interval, $L(a) = u$, and $L(b) = v$.  When the choice of geodesic arc
connecting $u$ to $v$ does not matter (or is clear from context) we
denote it by $[u, v]$.

This paper only considers {\em geodesic metric spaces}: metric spaces
in which every pair of points is connected by some geodesic arc.
However, we do not assume that our spaces are {\em proper}.

A subset $U \subset X$ is {\em quasi-convex} with constant $R \geq 0$
if, for every pair of points $u, v \in U$ and for every geodesic arc
$L$ connecting $u$ to $v$, the image of $L$ lies inside a closed $R$
neighborhood of the set $U$.

An arc $L \from [a, b] \to X$ is {\em quasi-geodesic} with constants
$\lambda \geq 1, \epsilon \geq 0$ if 
\begin{eqnarray}
\label{Eqn:QuasiGeodesic}
\frac{1}{\lambda} |b' - a'| - \epsilon \leq 
        d_X(L(a'), L(b')) \leq \lambda |b' - a'| + \epsilon
\end{eqnarray}
for every finite interval $[a', b'] \subset [a, b]$.

\subsection{Isometries}
A map $h \from X \to X$ is an {\em isometry} if for every pair of points
$x, y \in X$ we have $d_X(x, y) = d_X(h(x), h(y))$.  

\begin{define}
\label{Def:AveDisp}
Fix $x \in X$.  The {\em average displacement} of an isometry $h$ is
the quantity
$$ \avedisplace(h) = \lim_{n \to \infty} \frac{d_X(x, h^n(x))}{n} 
                   = \inf_{n \in \NN} \frac{d_X(x, h^n(x))}{n}. $$
\end{define}

It is well-known (see~\cite{CDP90}, Chapter~10, for example) that
$\avedisplace(h)$ exists and is independent of the given choice of
$x \in X$.  Note also that $n \cdot \avedisplace(h) \leq d_X(x, h^n(x))$, 
for all $n$.  We say an isometry is {\em hyperbolic} if
its average displacement is strictly positive.

\begin{remark}
This is one of several equivalent definitions of a hyperbolic isometry.
Claim~\ref{Clm:LOneIsQuasiGeodesic} below shows that any orbit of a
hyperbolic isometry, acting on a Gromov hyperbolic space, is a 
quasi-isometric embeddings of $\ZZ$, as expected.
\end{remark}

\subsection{Gromov hyperbolicity}
A geodesic metric space $(X, d_X)$ is {\em Gromov hyperbolic} with
constant $\delta$, or simply {\em $\delta$--hyperbolic}, if every
geodesic triangle is {\em $\delta$--thin}: the (closed) $\delta$
neighborhood of any two of the sides of the triangle contains the
third side.  Here a geodesic triangle is a collection of three
geodesic arcs which connect in pairs some triple of points $x, y, z
\in X$.  As an immediate corollary geodesic $n$--gons are $(n-2) \cdot
\delta$--thin: any one side lies in a $(n-2) \cdot \delta$ neighborhood
of the union of the other $n - 1$ sides.

Again and again we will use the remarkable ``stability'' property of
quasi-geodesics in a $\delta$--hyperbolic space:

\begin{lemma}
\label{Lem:SuckItIn}
For any choice of $\delta, \lambda, \epsilon$ there is a constant $R >
0$ such that: if $X$ is a $\delta$--hyperbolic space $X$ and $L \from
[a, b] \to X$ is a $(\lambda, \epsilon)$ quasi-geodesic, then the
image of $L$ is quasi-convex with constant $R$.  Furthermore, if $[a,
b]$ is a finite interval then the image of $L$ lies within the closed
$R$ neighborhood of any geodesic connecting $L(a)$ and $L(b)$.
\end{lemma}

See~\cite{Bridson99}, page~404, for a proof and note that $X$ need not
be proper.  The number $R$ is referred to as the {\em stability
constant} for $L$.

\section{Triangles and quadrilaterals}

Following~\cite{Bridson99}, page~463, we define a notion of a
``quasi-projection'' onto a quasi-convex set and deduce a few
consequences.

\subsection{Closest point projections}

Suppose that $U \subset X$ is nonempty and $X$ is $\delta$--hyperbolic.
Define a {\em quasi-projection} from $X$ to $U$ as follows: given
$\epsilon > 0$ and $y \in X$ put
$$\proj_U^\epsilon(y) = \{ y' \in U \st 
d_X(y, y') \leq d_X(y, U) + \epsilon \}.$$
That is, $\proj_U^\epsilon(y)$ is the set of points in $U$ which are,
within an error of $\epsilon$, closest to $y$.  Note that
$\proj_U^\epsilon(y)$ is nonempty.  


\begin{remark}
When $U$ is quasi-convex the function $\proj_U^\epsilon$ is a {\em
quasi-map} from $X$ to $U$.  That is, the diameter of
$\proj_U^\epsilon(y)$ is bounded independently of the point $y$.  This
is a direct consequence of Lemma~\ref{Lem:Triangles}, below.
\end{remark}

\subsection{The geometry of projections}

This section discusses similarities between the function
$\proj_U^\epsilon$ and orthogonal projection in hyperbolic space
$\HH^n$.

\begin{lemma}[Triangle lemma]
\label{Lem:Triangles}
Fix $\epsilon>0$.  Suppose $X$ is $\delta$--hyperbolic and $U \subset
X$ is quasi-convex with constant $R$.  Suppose that $y \in X$, that
$y' \in \proj_U^\epsilon(y)$, and that $u$ is another point of $U$.
Then
$$d_X(y, y') + d_X(y', u) \leq d_X(y, u) + (2\epsilon + 4\delta + 2R).$$ 
\end{lemma}

\begin{proof}
Note that if $d_X(y', u) \leq \epsilon + 2\delta + R$ then the
conclusion follows from the triangle inequality.  Suppose, then, that
$\epsilon + 2\delta + R < d_X(y', u)$ and let $a$ be the point of
$[y', u]$ such that $d_X(y', a) = \epsilon + 2\delta + R + \epsilon'$,
where $0 < \epsilon' < d_X(y', u) - \epsilon - 2\delta - R$.

\begin{claim}
\label{Clm:FarAway}
The point $a$ does not lie within a $\delta$ neighborhood of $[y', y]$.
\end{claim}

\begin{proof}[Proof of Claim~\ref{Clm:FarAway}]
Suppose the opposite.  Then there is a point $b \in [y', y]$ with
$d_X(a, b) \leq \delta$.  See Figure~\ref{Fig:Triangles}.  By the
triangle inequality we have $ \epsilon + \delta + R + \epsilon' \leq
d_X(y', b) $.  On the other hand consider a piecewise geodesic from
$y$ to $b$ to $a$ to $U$.  This has length at most $d_X(y, b) + \delta
+ R$ and at least $d_X(y', y) - \epsilon$.  Thus $d_X(y', y) -
\epsilon \leq d_X(y, b) + \delta + R$.

\begin{figure}[ht!]\small
\psfrag{x}{$y$}
\psfrag{x'}{$y'$}
\psfrag{u}{$u$}
\psfrag{a}{$a$}
\psfrag{b}{$b$}
\psfrag{c}{$c$}
$$\begin{array}{c}
\epsfig{file=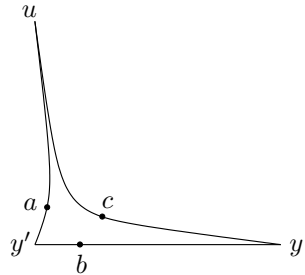, height = 3.5 cm}
\end{array}$$
\caption{A ``right-angled'' triangle}
\label{Fig:Triangles}
\end{figure}

\begin{remark}
This ``no short-cuts'' principle underlies both arguments in this
section.
\end{remark}

Subtract $d_X(y, b)$ from both sides of the above inequality to find
that $d_X(y', b)$ $\leq \epsilon + \delta + R$.  This contradiction proves
the claim.
\end{proof}

Continuing with the proof of Lemma~\ref{Lem:Triangles}, by 
$\delta$--thinness of the triangle $uy'y$ we have a point $c \in [y,
u]$ with $d_X(a, c) \leq \delta$.  See Figure~\ref{Fig:Triangles}.
Attempting to shortcut from $y$ to $c$ to $a$ to $U$ shows that
$d_X(y, y') \leq d_X(y, c) + \epsilon + \delta + R$.  Also, by the
triangle inequality applied to $uy'c$, we have $d_X(y', u) \leq d_X(u,
c) + \delta + \epsilon + 2\delta + R + \epsilon'$.  Adding these last
two inequalities and letting $\epsilon'$ tend to zero proves the
lemma.
\end{proof}

A similar lemma holds for quadrilaterals.

\begin{lemma}[Quadrilateral lemma]
\label{Lem:Quads}
Fix $\epsilon > 0$.  Suppose $X$ is a $\delta$--hyperbolic space and
that $U \subset X$ is quasi-convex with constant $R$.  Suppose that
$y, z \in X$ while $y' \in \proj_U^\epsilon(y)$ and $z' \in
\proj_U^\epsilon(z)$.  Suppose that $2\epsilon + 8\delta + 2R <
d_X(y', z')$.  Then
$$d_X(y, y') + d_X(y', z') + d_X(z', z)
            \leq d_X(y, z) + (4\epsilon + 12\delta + 4R).$$ 
\end{lemma}




See also Chapter 10, Proposition~2.1, of~\cite{CDP90}, Proposition
$\rm{III}.\Gamma.3.11$ of \cite{Bridson99}, or Lemma~7.3.D
of~\cite{Gromov87}.

\begin{proof}[Proof of Lemma~\ref{Lem:Quads}]
Let $a \in [y',z']$ be the point with $d_X(a, z') = \epsilon + 4\delta
+ R + \epsilon'$, where $0<\epsilon'< d_X(y',z')-2\epsilon-8\delta-2R$.  
As in Claim~\ref{Clm:FarAway}
the point $a$ lies outside of a $2\delta$ neighborhood of the union
of $[y, y']$ and $[z, z']$.  Thus, by $2\delta$--thinness of the
geodesic quadrilateral $yy'z'z$, there is a point $c \in [y, z]$ such
that $d_X(a, c) \leq 2\delta$.  See Figure~\ref{Fig:Quads}.

\begin{figure}[ht!]\small
\psfrag{x}{$z$}
\psfrag{x'}{$z'$}
\psfrag{y}{$y$}
\psfrag{y'}{$y'$}
\psfrag{a}{$a$}
\psfrag{c}{$c$}
\psfrag{c'}{$d$}
\psfrag{e}{$f$}
$$\begin{array}{c}
\epsfig{file=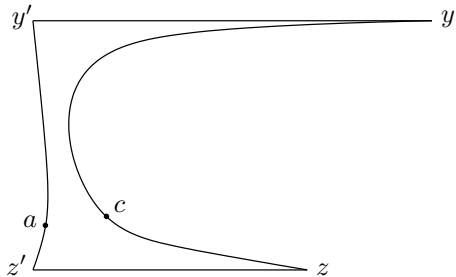, height = 3.5 cm}
\end{array}$$
\caption{Quadrilateral with two ``right angles''}
\label{Fig:Quads}
\end{figure}

Applying the triangle inequality to the piecewise geodesic from $y$ to
$c$ to $a$ to $z'$ we find that $d_X(y, z') \leq d_X(y, c) + 2\delta +
\epsilon + 4\delta + R + \epsilon'$.  Combining this with
Lemma~\ref{Lem:Triangles} (taking $u = z'$) gives $d_X(y, y') +
d_X(y', z') \leq d_X(y, c) + 3\epsilon + 10\delta + 3R + \epsilon'$.  Finally,
short-cutting from $z$ to $c$ to $a$ to $U$ gives $d_X(z, z') \leq d_X(z,
c) + \epsilon + 2\delta + R$.  Adding these last two inequalities
and letting $\epsilon'$ tend to zero gives the desired conclusion.
\end{proof}

\section{The displacement lemma}
\label{Sec:Displacement}

This section gives a proof of:

\begin{lemma}[Displacement lemma]
\label{Lem:Displacement}
Suppose $h$ is a hyperbolic isometry of the $\delta$--hyperbolic space
$X$.  Suppose that $x \in X$ is fixed.  Then there is a constant $K$
such that, for all $n$,
$$ | d_X(x, h^n(x)) - n \cdot \avedisplace(h) | \leq K. $$
\end{lemma}

Recall that $\avedisplace(h)$ is the average displacement of $h$; see
Definition~\ref{Def:AveDisp}.  Though this result seems to be
well-known, we have not been able to find it in the literature.

\begin{remark}
Lemma~\ref{Lem:Displacement} is immediate in $\HH^n$.  In fact
$\delta$--hyperbolicity is not necessary; the displacement lemma also
holds for semi-simple isometries of convex metric spaces.
See~\cite{Bridson99} for definitions.
\end{remark}

\begin{remark}
Set $\displace_n = d_X(x, h^n(x))/n$.  A consequence of the lemma is
that $\displace_n - \avedisplace(h) = O(1/n)$.
\end{remark}

Throughout the proof of the lemma we adopt the notation $x_n = h^n(x)$
where $x_0 = x$ is the basepoint provided by the hypothesis.  Note
that $n\cdot \avedisplace(h) - d_X(x_0, x_n) \leq 0$.  This gives the
upper bound.  For the lower, we define a sequence of infinite arcs in
$X$: for each positive integer $n$ fix a geodesic arc
$$P_n = [x_0, x_n].$$
Let
\begin{equation}
\label{Eqn:Ln}
L_n = \bigcup_{k \in \ZZ} h^{nk}(P_n)
\end{equation}
be parameterized by arc-length.  We will investigate these infinite
arcs in order to prove the lemma.

Again, take $\displace_m = d_X(x_0, x_m)/m$.  Fix $M_0 \geq 3$ so that
if $m + 1 \geq M_0$ then $|\displace_m - \avedisplace(h)| <
\frac{1}{2}\avedisplace(h)$.  (See Definition~\ref{Def:AveDisp}.)

\begin{claim}
\label{Clm:LOneIsQuasiGeodesic}
The arc $L_1$ is a quasi-geodesic with constants $\lambda_1 =
\frac{2\displace_1}{\avedisplace(h)}$ and $\epsilon_1 = M_0
\displace_1$.
\end{claim}

See also~\cite{CDP90}, Chapter~10, Proposition~6.3.

\begin{remark}
More is true.  All of the arcs $L_n$ are quasi-geodesic
with uniformly bounded additive constants $\epsilon_n$ and
multiplicative constants of the form $\lambda_n = 1 + O(1/n)$.
\end{remark}

\begin{proof}[Proof of Claim~\ref{Clm:LOneIsQuasiGeodesic}]
Fix $a < b \in \RR$.  We will show that $d_X(L_1(a), L_1(b))$
satisfies the inequalities given by Equation~\ref{Eqn:QuasiGeodesic}.
To simplify notation set $L_1(a) = u$ and $L_1(b) = v$.  As $L_1$ is
parameterized by arc-length, $d_X(u, v) \leq |b - a|$.

By applying $h$ some number of times we may assume that $u$ lies in
$P_1$ while $v$ lies in $h^m(P_1)$, for some smallest possible
nonnegative integer $m$.  Recall that $x_m = h^m(x)$.  The triangle
inequality gives
$$d_X(x_0, x_m) \leq d_X(u, v) + 2\displace_1$$
where $\displace_1 = d_X(x_0, x_1)$.  In the case $m + 1 \geq M_0$ we have
$$\frac{m}{2} \avedisplace(h) < d_X(x_0, x_m)$$
and, as $L_1$ is parameterized via arc-length, 
$$|b - a| \leq (m + 1) \displace_1.$$
Chaining together the above three inequalities gives:
$$\frac{\avedisplace(h)}{2\displace_1}|b - a| - \frac{\avedisplace(h)}{2} 
                                - 2\displace_1 \leq d_X(u, v).$$

On the other hand, if $m + 1 \leq M_0$ then
$$|b - a| - M_0 \displace_1 \leq d_X(u, v).$$
Since $\avedisplace(h) \leq \displace_1$ and $M_0 \geq 3$, regardless
of $m$ we have
$$\frac{\avedisplace(h)}{2\displace_1} |b - a| - M_0 \displace_1 \leq
d_X(u, v).$$
This completes the proof of the claim. 
\end{proof}

\begin{define}
\label{Def:Axis}
We call $L_1$ as defined in Equation~\ref{Eqn:Ln} an {\em invariant
axis} for $h$.
\end{define}

Returning to the proof of Lemma~\ref{Lem:Displacement}, choose $n \in
\NN$ and fix attention on the infinite arc $L_n$.  Recall that $x_m =
h^m(x)$ where $x$ is the chosen basepoint.  As $\displace_m = d_X(x_0,
x_m)/m$ converges to $\avedisplace(h)$ from above there is a positive
integer $M_1$ such that if $m > M_1$ then $\displace_m -
\avedisplace(h) \leq 1/n$.  Choose $m > M_1$ of the form $m = kn$.  It
follows that
\begin{eqnarray}
\label{Eqn:MEqualsKN}
\frac{d_X(x_0, x_m)}{k} - n \cdot \avedisplace(h) \leq 1.
\end{eqnarray}

\begin{figure}[ht!]\small
\psfrag{x}{$x_0$}
\psfrag{hnx}{$x_n$}
\psfrag{hmx}{$x_m$}
\psfrag{h2nx}{$x_{2n}$}
\psfrag{z1}{$z_1$}
\psfrag{z2}{$z_2$}
$$\begin{array}{c}
\epsfig{file=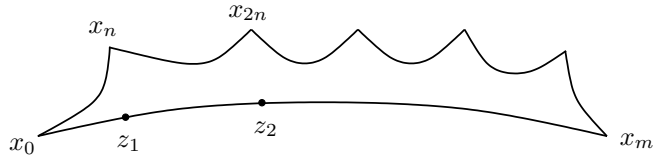, height = 2 cm}
\end{array}$$
\caption{The sides of the ladder:  recall that $x_k = h^k(x_0)$.  The
  rungs connecting $x_{ni}$ to $z_i$ are not drawn.}
\label{Fig:LadderClimb}
\end{figure}

Now we compare the quantities $d_X(x_0, x_m)/k$ and $d_X(x_0, x_n)$.
Recall that $P_m = [x_0, x_m]$.  For $i \in \{0, 1, \ldots k\}$
choose a point $z_i \in \proj_{P_m}^0(x_{ni})$.  Note that $z_i$
exists as $P_m$ is compact.
It follows that $d_X(x_{ni}, z_i) \leq R$ where $R$ is the
stability constant provided by Lemma~\ref{Lem:SuckItIn} 
for the quasi-geodesic $L_1$.  See Figure~\ref{Fig:LadderClimb}.

\begin{claim}
\label{Clm:NoFolding}
If $n > \frac{3R}{\avedisplace(h)}$ then, when traveling from $z_0 =
x_0$ to $z_{k} = x_m = h^m(x)$ along $P_m$, the points $z_i$ are
distinct and encountered in order of their indices.
\end{claim}

\begin{proof}[Proof of Claim~\ref{Clm:NoFolding}]
From the definition of $\avedisplace(h)$ we know that $d_X(x_0,
x_{nj}) \geq n \cdot \avedisplace(h) > 3jR$.  Recall that $d_X(z_i,
x_{ni}) \leq R$.  So for all $i$, $d_X(z_i, z_{i+j}) \geq 3jR - 2R$.
Thus consecutive $z_i$'s are distinct.

If the $z_i$'s are not encountered in order by index then there is a
smallest $i$ with $z_{i+1}$ appearing before $z_i$.  So either
$z_{i+1} \in [z_{i-1}, z_i]$ or $z_{i-1} \in [z_{i+1}, z_i]$.  Suppose
that former occurs (the latter is similar).  Now, by the above
paragraph, $d_X(z_{i-1}, z_{i}) \leq d_X(x_0, x_n) + 2R$ and
$d_X(z_{i+1}, z_{i}) \geq d_X(x_0, x_n) - 2R$.  We deduce that
$d_X(z_{i-1}, z_{i+1}) \leq 4R$.  Thus by the triangle inequality we
have $d_X(x_{n(i-1)}, x_{n(i+1)}) \leq 6R$.  But this is a
contradiction.
\end{proof}

To conclude the proof of Lemma~\ref{Lem:Displacement}, assume
for the moment that $n > \frac{3R}{\avedisplace(h)}$.  By the
above claim the intervals $[z_i, z_{i+1}]$ form a disjoint partition
of $P_m$ (ignoring endpoints).  Let $j \in \{ 0, 1, \ldots, k-1 \}$ be
chosen to minimize the length of $[z_j, z_{j+1}]$.  Then certainly
$d_X(z_j, z_{j+1}) \leq d_X(x_0, x_m)/k$.  Applying the triangle
inequality to the rectangle with vertices $x_{nj}$, $z_{j}$,
$z_{j+1}$, $x_{n(j+1)}$ shows that
$$d_X(x_0, x_n) \leq \left(d_X(x_0, x_m)/k\right) + 2R.$$
Conclude, by Equation~\ref{Eqn:MEqualsKN}, that 
$$d_X(x_0, x_n) - n \cdot \avedisplace(h) \leq 1 + 2R$$
as long as $n > \frac{3R}{\avedisplace(h)}$.  Thus regardless of $n$
$$|d_X(x_0, x_n) - n \cdot \avedisplace(h)| \leq \max\left\{1 + 2R,
\frac{3R\displace_1}{\avedisplace(h)} \right\}.$$ 
This proves the Displacement Lemma. \qed

\section{Bounded projection implies linear growth}
\label{Sec:BoundedProj}

Let $X$ be a $\delta$--hyperbolic space and $h$ a hyperbolic isometry
of $X$.  Choose an invariant axis $L_1$ for $h$ as in
Definition~\ref{Def:Axis} above.

\begin{define}
\label{Def:BoundedProj}
A subset $Y \subset X$ has {\em bounded projection} with
respect to $h$ if there is $\epsilon > 0$
so that the set $P = \bigcup_{z \in Y}
\proj_{L_1}^\epsilon(z)$ has finite diameter.  
\end{define}

Note that (for given $Y$ and $h$) the above definition is independent
of the choices involved in defining $L_1$.  In this section we prove:

\begin{theorem}
\label{Thm:BoundedImpliesLinear}
If $Y$ has bounded projection with respect to $h$ then there is
a constant $K$ such that:
$$ | d_X(Y, h^n(Y)) - n \cdot \avedisplace(h) | \leq K. $$
\end{theorem}

\begin{proof}
Pick any point $y \in Y$.  Then there is a constant $K'$,
provided by Lemma~\ref{Lem:Displacement}, such that $d_X(y, h^n(y)) -
n \cdot \avedisplace(h) \leq K'$ for all $n$.  So $d_X(Y,
h^n(Y)) - n \cdot \avedisplace(h) \leq K'$ for all $n$.
This gives the desired upper bound.  We now turn to the lower bound.

Let $y \in Y$ and take $z \in h^n(Y)$.  We must bound from below the
quantity $d_X(y, z)$.  Choose $y' \in \proj_{L_1}^\epsilon(y)$ and $z'
\in \proj_{L_1}^\epsilon(z)$.  Let $P = \bigcup_{w \in Y}
\proj_{L_1}^\epsilon(w)$ be the projection of $Y$ to $L_1$.  By
hypothesis this set has finite diameter, say $K''$.  Note
that $d_X(z', h^n(y')) \leq K''$.

By Definition~\ref{Def:AveDisp} we have $n \cdot \avedisplace(h) \leq
d_X(y', h^n(y'))$.  Thus, by the triangle inequality, $n \cdot
\avedisplace(h) \leq d_X(y', z') + K''$.

For the remainder of the proof take $M \equiv \frac{2\epsilon +
  8\delta + 2R + K''}{\avedisplace(h)}$.  Then, when $n > M$ the
previous inequality implies that $2\epsilon + 8\delta + 2R < d_X(y',
z')$.  Letting $K''' = 4\epsilon + 12\delta + 4R$ and applying
Lemma~\ref{Lem:Quads} we have
\begin{eqnarray*}
d_X(y, z) + K''' 
          & \geq & d_X(y, y') + d_X(y', z') + d_X(z', z) \\
          & \geq & d_X(y', z').
\end{eqnarray*}


Thus when $n > M$, regardless of our choice of $y \in Y$ and $z \in
h^n(Y)$, we have
$$n \cdot \avedisplace(h) \leq d_X(y, z) + K'' + K'''.$$
That is,
$$n \cdot \avedisplace(h) - d_X(Y, h^n(Y)) \leq K'' + K'''$$ 
for $n > M$.  
Also, as discussed in the first paragraph of the proof,
$$d_X(Y, h^n(Y)) - n \cdot \avedisplace(h) \leq K'$$
for all $n$.  So take 
$$K = \max
    \left\{ K', K'' + K''', |d_X(Y,
    h^n(Y)) - n \cdot \avedisplace(h)| \right\}$$ 
where $n$ ranges from $1$ to $M$. 
This gives the desired bound.
\end{proof}

\section{B implies Lin}
\label{Sec:BImpliesLin}

This section transforms the preceding purely geometric
considerations into tools appropriate to the setting of Heegaard
splittings.  We begin by defining the graph of curves.

Let $S$ be a closed orientable surface of genus at least two.  Let
$\mathcal{C}^0(S)$ be the set of isotopy classes of essential simple
closed curves in $S$.  The {\em graph of curves}, $\mathcal{C}^1(S)$,
has vertex set $\mathcal{C}^0(S)$ and an edge connecting two distinct
vertices if and only if the two curves may be realized disjointly.  We
take each such edge to be a copy of the interval $[0, 1]$ and give
$\mathcal{C}^1(S)$ the induced metric.

\begin{remark}
The graph $\mathcal{C}^1(S)$ is the one-skeleton of the {\em curve
complex} and is quasi-isometric to the full complex.  For
simplicity we consider only $\mathcal{C}^1(S)$.
\end{remark}

We require a pair of deep results of Masur and
Minsky~\cite{MasurMinsky99}.

\begin{theorem}[Masur--Minsky]
\label{Thm:CIsDeltaHyperbolic}
The graph of curves $(\mathcal{C}^1(S), d_\mathcal{C})$, is a Gromov
hyperbolic space.  Furthermore, pseudo-Anosov maps act on the graph of
curves as hyperbolic isometries.
\end{theorem}

Thus, Theorem~\ref{Thm:BoundedImpliesLinear} can be translated to the
language of Heegaard splittings as follows.  Recall that $\mathcal{V}
\subset \mathcal{C}^0(S)$ is the set of curves which bound essential
disks in the handlebody $V$.  As usual $S = \bdy V$ has genus two
or more.  Fix $h \from S \to S$ a pseudo-Anosov map.

\begin{define}
\label{Def:B}
Given $V$, $S$, and $h$ as above, we say that \Bounded\ holds if the
handlebody set $\mathcal{V}\subset \mathcal{C}^1(S)$ has bounded
projection with respect to $h$.
\end{define}

Also, note that the distance of a Heegaard splitting $(S, \mathcal{V},
\mathcal{W})$ (Definition~\ref{Def:Distance}) is exactly the quantity
$d_\mathcal{C}(\mathcal{V}, \mathcal{W})$.  

\begin{define}
\label{Def:Lin}
Given $V$, $S$, and $h$ as above, we say that \Lin\ holds if there is
a constant $K$ so that, for all $n \in \NN$,
$$|d_\mathcal{C}(\mathcal{V}, h^n(\mathcal{V})) 
             - n \cdot \avedisplace(h)| \leq K.$$
\end{define}

\begin{corollary}
\label{Cor:BImpliesLin}
Suppose a handlebody $V$ with $S = \bdy V$ and a pseudo-Anosov map
$h \from S \to S$ are given.  Then \Bounded\ implies \Lin.
\end{corollary}

\begin{proof}
This follows immediately from Theorems~\ref{Thm:BoundedImpliesLinear}
and~\ref{Thm:CIsDeltaHyperbolic}.
\end{proof}

\begin{remark}
\label{Rem:Estimation}
Algorithmic computation of $d_\mathcal{C}(\mathcal{V},
h^n(\mathcal{V}))$ would be highly desirable.  Upper and lower bounds
for $\avedisplace(h)$ may perhaps be obtained using methods similar
to~\cite{MasurMinsky99}\footnote{Added in proof:
K~Shackleton~\cite{Shackleton04} using work of
B~Bowditch~\cite{Bowditch03} has made significant progress on this
problem.}.  Estimation of $K$ seems more difficult.  For any $n \in
\NN$ there are pairs $(V, h)$ where the projection of $\mathcal{V}$ to
$L_1$ has finite diameter which is bigger than $n$.


\end{remark}

\section{Equivalence of H and B}
\label{Sec:HIffB}

In this section we deduce the following from work of
Klarreich~\cite{Klarreich99}.

\begin{theorem}
\label{Thm:HIffB}
Suppose a handlebody $V$ with $S = \bdy V$ and a pseudo-Anosov map $h
\from S \to S$ are given.  Then \Hempel\ is equivalent to \Bounded.
\end{theorem}

\subsection{The Gromov product}
Before proving Theorem~\ref{Thm:HIffB} we recall the definition of the
Gromov product.  Suppose that a basepoint $x_0$ in the
$\delta$--hyperbolic space $X$ is given.  The {\em Gromov product} of a
pair of points $y, z \in X$ is the quantity
$$( y \cdot z ) = \frac{1}{2} \left( d_X(x_0, y) + d_X(x_0, z) -
  d_X(y, z) \right).$$
Following~\cite{Gromov87} we say that a sequence $\{y_i\}_{i =
0}^{\infty} \subset X$ {\em converges at infinity} if $\lim_{n,m \to
\infty} ( y_n \cdot y_m )$ is infinite.
This is independent of the choice of basepoint, $x_0$.  Two such
sequences $\{y_i \}$ and $\{ z_i \}$, both converging at infinity, are
{\em equivalent} if $\lim_{n,m \to \infty} ( y_n \cdot z_m )$ is again
infinite.  The {\em Gromov boundary} of $X$, denoted $\bdy_{\infty}
X$, is the set of equivalence classes of sequences which converge at
infinity.  As a final bit of notation, set $|y| = (y \cdot y) =
d_X(x_0, y)$.

Now, if $h$ is a hyperbolic isometry of $X$ we define the {\em stable
  and unstable fixed points} of $h$ to be the points of $\bdy_\infty
X$ containing the sequences $L^+(h) = \{x_n = h^n(x_0) \st n \in
\NN\}$ and $L^-(h) = \{x_{-n} = h^{-n}(x_0) \st n \in \NN\}$
respectively.  Recall that $L_1$ is a piecewise geodesic through the
points $\{x_i\}_{i = -\infty}^\infty$.

There is a simple relation between projection to the quasi-geodesic
$L_1$ and the Gromov product.  Fix $\epsilon > d_X(x_0, x_1)$.  Again
we use $x_0$ as the basepoint for computing the Gromov product.

\begin{lemma}
\label{Lem:ProjToLOne}
Fix $y, z \in X$.  Pick $x_m \in \proj_{L_1}^{\epsilon}(y)$, $x_n \in
\proj_{L_1}^{\epsilon}(z)$, and suppose that $n, m > 0$.  Then
$$(y \cdot z) \geq \min\{ |x_n|, |x_m| \} - K$$  and
$$\left| (y \cdot x_n) - \min\{ |x_n|, |x_m| \} \right| \leq K.$$
Here $K$ is a constant
not depending on $y$ or $z$.
\end{lemma}

\begin{proof}
Suppose that $m \leq n$ as the other case is similar.  For ease of
notation let $A = d_X(x_0, x_m) = |x_m|$, $B = d_X(x_m, y)$, $C =
d_X(x_m, x_n)$, $D = d_X(x_n, z)$, and $E = d_X(x_0, x_n) = |x_n|$.
Let $K' = 2\epsilon + 4\delta + 2R$, where $R$ is the stability
constant for $L_1$.

Now $(y \cdot z) = \frac{1}{2}(d_X(x_0, y) + d_X(x_0, z) - d_X(y,
z))$.  The first term is greater than $A + B - K'$, applying
Lemma~\ref{Lem:Triangles} and the fact that $L_1$ is quasi-geodesic
(Claim~\ref{Clm:LOneIsQuasiGeodesic}), hence quasi-convex.  Similarly,
the second term is greater than $D + E - K'$.  But $E \geq C + A -
2R$, using the triangle inequality and the fact that $L_1$ is a
quasi-geodesic.  Finally, the third term is less than $B + C + D$
using the triangle inequality.  So $(y \cdot z) \geq \frac{1}{2}(A + B
- K' + A + C + D - 2R - K' - B - C - D)$ and we have the desired lower
bound.

When $z = x_n$ we also obtain an upper bound, as in this case the
first term is less than $A + B$, the second is less than $A + C$,
while the third is greater than $B + C - K'$.
\end{proof}
This leads to:
\begin{lemma}
\label{Lem:NotBIffConverge}
A subset $Y \subset X$ has unbounded projection with respect to
$h$ if and only
if there is a sequence $y_n \in Y$ converging to the stable or
unstable fixed point for $h$ at infinity.
\end{lemma}

\begin{proof}
Suppose that the set $Y$ has unbounded projection to the sequence
$L^+(h) = \{x_n = h^n(x_0) \st n \in \NN \}$.  (The other case is
similar.)  Choose $\epsilon > d_X(x_0, x_1)$ and a sequence $\{y_n\}
\subset Y$ so that $x_{m(n)} \in \proj_{L_1}^\epsilon(y_n)$, for some
$m(n) > n$.  It follows from Lemma~\ref{Lem:ProjToLOne} that $(y_k
\cdot y_l) \geq \min\{ |x_{m(k)}|, |x_{m(l)}|\} - K$ and thus
$\{y_n\}$ converges at infinity.  Also, using the second inequality of
Lemma~\ref{Lem:ProjToLOne} it is easy to check that $\{y_n\}$ and
$\{x_{m(n)}\}$ are equivalent.

On the other hand, suppose that there is a sequence $\{y_n\} \subset
Y$ with $\{y_n\}$ converging to the stable fixed point of $h$.  Then
$\{y_n\}$ and $\{x_n\}$ are equivalent.  So we may pass to
subsequences $\{y_k\}$ and $\{x_k\}$ so that $(y_k \cdot x_k)$ goes to
infinity with $k$.  Pick $x_{m(k)} \in \proj_{L_1}^\epsilon(y_k)$.
Then, by the second half of Lemma~\ref{Lem:ProjToLOne}, the quantity
$\min \{|x_{m(k)}|, |x_k|\}$ must also tend to infinity with $k$.
Thus $|x_{m(k)}|$ also tends to infinity with $k$ and we are done.
\end{proof}

\subsection{The boundary of the curve complex}
We next cite the necessary component from
Klarreich~\cite{Klarreich99}.  Let $\MinLam$ be the space of minimal
measured laminations, considered up to topological equivalence (ie,
take a quotient by forgetting the measures).  Klarreich gives a
homeomorphism $\pi \from \MinLam \to \bdy_\infty \mathcal{C}^1(S)$
such that the following holds\footnote{Added in proof:
U~Hamenstaedt~\cite{Hamenstaedt04} has announced a combinatorial
proof of Klarreich's theorem.}:

\begin{theorem}[Klarreich]
\label{Thm:Klarreich}
Let $\gamma_n$ be a sequence of essential simple closed curves in the
surface $S$.  Suppose that $\mathcal{L}$ is a minimal, uniquely
ergodic lamination on $S$.  Then the sequence $\gamma_n \in \PML(S)$
converges to $\mathcal{L}$ if and only if $\gamma_n \in
\mathcal{C}^1(S)$ converges to $\pi(\mathcal{L}) \in \bdy_\infty
\mathcal{C}^1(S)$.
\end{theorem}

See Theorem~3.2 of~\cite{Minsky01}, for a more precise version of
Klarreich's result.  We are now ready to prove
Theorem~\ref{Thm:HIffB}.

\begin{proof}[Proof of Theorem~\ref{Thm:HIffB}]
To begin, pick any $x \in \mathcal{C}^0(S)$, and let $L_1$ be a
quasi-geodesic (as defined in Section~\ref{Sec:Displacement}) passing
through the points $\{ x_n = h^n(x) \st n \in \ZZ \}$.  Also take $x =
x_0$ to be the basepoint when computing the Gromov product.

Note that \Hempel\ does not hold if and only if there is a sequence of
curves $v_n \in \mathcal{V}$ such that $v_n$ converges in $\PML(S)$ to
one of $\mathcal{L}^{\pm}(h)$.  Suppose $v_n$ converges to
$\mathcal{L}^+(h)$. (The other case is identical.)

Applying Theorem~\ref{Thm:Klarreich} the $v_n$ converge in $\PML(S)$ to
$\mathcal{L}^+(h)$ if and only if they also converge, in
$\mathcal{C}^1(S)$, to $\pi(\mathcal{L}^+(h))$ in the boundary of the
curve complex.  Now apply Lemma~\ref{Lem:NotBIffConverge} to find that
this occurs if and only if $\mathcal{V}$ has unbounded projection to
$L^+(h)$.
\end{proof}

\section{The converse: Unb implies B}
\label{Sec:UnbImpliesB}

In this section we prove Theorem~\ref{Thm:UnbImpliesB}, \Unb\ implies
\Bounded.  Recall that $\mathcal{V} \subset \mathcal{C}^0(S)$, the
handlebody set, contains all curves which bound disks in the
handlebody $V$.  We will need a final result of Masur and
Minsky~\cite{MasurMinsky03}:

\begin{theorem}[Masur--Minsky]
\label{Thm:VIsQuasiConvex}
Fix a handlebody $V$ with $\bdy V=S$.  The handlebody set
$\mathcal{V}$ is quasi-convex in $\mathcal{C}^1(S)$.
\end{theorem}

We need one preparatory lemma about quasi-convex sets.

\begin{lemma}
\label{Lem:YZMeetImpliesClose}
Suppose $X$ is a $\delta$--hyperbolic space and $Y$ and $Z$ are
quasi-convex subsets with constant $R$.  There is a constant $K$,
depending only on $\delta$ and $R$, such that: if $\{y_m\} \subset Y$
and $\{z_n\} \subset Z$ converge to the same point of $\bdy_\infty X$
then $d_X(Y, Z) < K$.
\end{lemma}

To paraphrase: if $Y$ and $Z$ are quasi-convex and intersect at
infinity then $Y$ and $Z$ are close to each other.

\begin{proof}
Pick $y_0$ to be the basepoint for computing the Gromov product.  Set
$D = d_X(y_0, z_0)$.  As $\lim_{m,n \to \infty} (y_m \cdot z_n) =
\infty$ we also have $\lim_{m \to \infty} (y_m \cdot z_m) =
\infty$.  Thus there is a large $k > 0$ so that any geodesic $[y_k, z_k]$
lies outside of a $5\delta + D$ ball about $y_0$.

Consider, then, a geodesic quadrilateral with vertices, in order,
$y_0$, $z_0$, $z_k$, and $y_k$.  Recall that the quadrilateral is
$2\delta$--thin: any one side lies in the $2\delta$ neighborhood of the
union of the other three sides.  Apply this to the side $[y_0, y_k]$.
But the $2\delta$ neighborhood of $[y_0, z_0]$ lies within the
$2\delta + D$ ball about $y_0$.  Also, the $2\delta$ neighborhood of
$[y_k, z_k]$ lies without the $3\delta + D$ ball about $y_0$.  Thus
some point of $[y_0, y_k]$ lies within $2\delta$ of some point of
$[z_0, z_k]$.  Apply quasi-convexity to find that $d_X(Y, Z) < 2\delta
+ 2R$.
\end{proof}

We may now prove:

\begin{theorem}
\label{Thm:UnbImpliesB}
Suppose a handlebody $V$ with $S = \bdy V$ and a pseudo-Anosov map $h
\from S \to S$ are given.  Then \Unb\ implies \Bounded.
\end{theorem}

\begin{proof}
We prove the contrapositive.  Choose $x \in \mathcal{C}^0(S)$ to be
the basepoint and construct $L_1$, a piecewise geodesic connecting the
points $\{x_n = h^n(x)\}$, as in Section~\ref{Sec:Displacement}.
Suppose that $\mathcal{V}$ has unbounded projection to $L_1$.  Thus,
by Lemma~\ref{Lem:NotBIffConverge}, there is a sequence $v_m \in
\mathcal{V}$ converging to the stable fixed point for $h$ at infinity.
(The case where $v_m$ converges to the unstable fixed point for $h$ is
identical.)

Note that, for any fixed $n \in \ZZ$, the same holds of the handlebody
set $h^n(\mathcal{V})$, as unbounded projection to $L_1$ is
$h$--invariant.  By Theorem~\ref{Thm:VIsQuasiConvex} both $\mathcal{V}$
and $h^n(\mathcal{V})$ are quasi-convex.  Thus both requirements of
Lemma~\ref{Lem:YZMeetImpliesClose} are satisfied.  It follows that the
distance $d_\mathcal{C}(\mathcal{V}, h^n(\mathcal{V}))$ is bounded
independently of $n$.
\end{proof}

\section{Equivalence of M and H}
\label{Sec:MIffH}

As usual we have a handlebody $V$ of genus at least two, $S = \bdy V$,
and $h$ a pseudo-Anosov map on $S$.  The following subset was
introduced by Masur in~\cite{Masur86}.

\begin{define}
\label{Def:MasurDomain}
The {\em Masur domain} $\mathcal{M}(V) \subset \PML(S)$ is the set of
laminations having nonzero geometric intersection with every
lamination in the closure of $\mathcal{V}$.
\end{define}

\begin{define}
\label{Def:M}
Given $V$, $S$, and $h$ as above, we say that \Masur\ holds if both
the stable and unstable laminations of $h$ lie in $\mathcal{M}(V)$.
\end{define}

We have:

\begin{lemma}
\label{Lem:MIffH}
Suppose a handlebody $V$ with $S = \bdy V$ and a pseudo-Anosov map $h
\from S \to S$ are given.  Then \Masur\ is equivalent to \Hempel.
\end{lemma}

\begin{proof}
Any lamination has zero geometric intersection with itself.  So the
closure of $\mathcal{V}$ lies in the complement of $\mathcal{M}(V)$.
Thus \Masur\ implies \Hempel.

Now consider the stable or unstable lamination of $h$,
$\mathcal{L}^\pm(h)$.  As $\mathcal{L}^\pm(h)$ is minimal and uniquely
ergodic, $\iota(\mu, \mathcal{L}^\pm(h)) = 0$ implies that $\mu =
\mathcal{L}^\pm(h)$ as projective measured laminations.  Thus if
$\mathcal{L}^\pm(h)$ is not in the Masur domain of $V$ then it must be
in the closure of $\mathcal{V}$.
\end{proof}

\section{Equivalence of K and H}
\label{Sec:KIffH}

In~\cite{Kobayashi88b} Kobayashi gives several examples of pairs
$(V,h)$ satisfying a condition which we call \Kobayashi\ (see
Definition~\ref{Def:K}).  In this section we prove that \Kobayashi\ is
equivalent to Hempel's hypothesis \Hempel.

\subsection{Laminations and the Whitehead graph}
A {\em cut system} $C$ for a handlebody $V$ is a collection of
disjoint, nonparallel, essential disks in $V$ so that the closure of
$V \setminus C$ (in the path metric) is a union of three-balls.  A
maximal cut system will be referred to as a {\em pants decomposition}
of $V$.  If $V \setminus C$ is connected then $C$ is a {\em minimal}
cut system.  Recall that $S = \bdy V$.  Let $\mathcal{L}$ be a
measured lamination on $S$.

\begin{define}
\label{Def:Tight}
The lamination $\mathcal{L}$ is {\em tight} with respect to a cut
system $C$ if no component of $S \setminus (C \cup \mathcal{L})$ is
{\em trivial:} has boundary a union of two arcs, $\alpha \cup \beta$,
where $\alpha \cap \beta = \bdy \alpha = \bdy \beta$, $\alpha \subset
\mathcal{L}$ and $\beta \subset \bdy C$.
\end{define}


When $C$ and $\mathcal{L}$ are tight we form the {\em Whitehead graph}
$\Gamma(\mathcal{L}, C)$ as follows.  Let $P$ be the closure of $S
\setminus \bdy C$, in the path metric.  Thus $P$ is a disjoint union
of planar surfaces.  Every arc of $\mathcal{L} \cap P$ now falls into
one of finitely many homotopy classes of properly embedded arcs in $P$.
The vertices of $\Gamma(\mathcal{L}, C)$ are the boundary components
of $P$.  For every homotopy class of arc we have an edge with the
obvious endpoints.  Note that every such arc inherits a positive
transverse measure from $\mathcal{L}$.  Note that the number of
components of $\Gamma(\mathcal{L}, C)$ equals the number of components
of $P$, when $\mathcal{L}$ is minimal.

A component of $\bdy P$ is a {\em cut vertex} for the graph if
removing the component, and all (open) edges adjacent to it, increases the
number of components of $\Gamma(\mathcal{L}, C)$.

\begin{define}
\label{Def:FullType}
A minimal lamination $\mathcal{L}$ is of {\em full type} with respect
to $C$ if $\mathcal{L}$ and $C$ are tight and the associated Whitehead
graph has no cut vertex.
\end{define}
\noindent 
One obstruction to being full type is the presence of {\em waves}.

\begin{define}
\label{Def:Wave}
A {\em wave} is any component of $\mathcal{L} \setminus C$ giving a
loop-edge in the Whitehead graph $\Gamma(\mathcal{L}, C)$.
\end{define}

\subsection{Kobayashi's hypothesis}
Let $V$ be a handlebody (of genus at least two) and $h$ a 
pseudo-Anosov map from $S = \bdy V$ to itself.  
From~\cite{Kobayashi88b} we have:

\begin{define}
\label{Def:K}
Given $V$, $S$, and $h$ as above, we say that \Kobayashi\ holds if
there are pants decompositions $C^\pm$ such that $\mathcal{L}^\pm(h)$
is of full type with respect to $C^\pm$.
\end{define}

\begin{theorem}
\label{Thm:KIffH}
Suppose a handlebody $V$ with $S = \bdy V$ and a pseudo-Anosov map $h
\from S \to S$ are given.  Then \Kobayashi\ is equivalent to \Hempel.
\end{theorem}
\noindent
We prove the two directions separately.

\begin{lemma}
\label{Lem:KImpliesH}
If \Kobayashi\ holds then \Hempel\ holds.
\end{lemma}
 

\begin{proof}
Suppose \Kobayashi\ holds.  We fix attention on $\mathcal{L}^+(h)$
and the maximal cut system $C^+$, as the other case is identical.
Let $v_n$ be a sequence from $\mathcal{V}$, converging in $\PML(S)$
to some minimal lamination, $\mathcal{L}$.  We must show that
$\mathcal{L} \neq \mathcal{L}^+(h)$.

Isotope each of the $v_n$'s so that each is tight with respect to
$C^+$.  Let $P$ be the collection of pants obtained by cutting $S$
along $\bdy C^+$ and taking the closure in the path metric.  Passing
to a subsequence, if necessary, we may assume that every pair $(v_n,
C^+)$ yields the same Whitehead graph, $\Gamma = \Gamma(v_n, C^+)$.
Note that $\Gamma$ contains $\Gamma' = \Gamma(\mathcal{L}, C^+)$ as a
subgraph.

As the $v_n$'s bound disks, there is some component $\rho$ of $\bdy P$ 
so that every $v_n$ contains a wave for $\rho$.  It follows that
$\rho$ is a cut vertex for $\Gamma$ and hence for $\Gamma'$.  So
$\mathcal{L}$ is not of full type and cannot be equal to
$\mathcal{L}^+(h)$.  
\end{proof}

The converse is more difficult and is dealt with in two steps.

\begin{lemma}
\label{Lem:HImpliesK'}
If \Hempel\ holds then there are minimal cut systems $C^\pm$ so that
$\mathcal{L}^\pm(h)$ is full type with respect to $C^\pm$.
\end{lemma}
 
\begin{proof}
We prove the contrapositive.  Suppose, as the other case is similar,
that the stable lamination $\mathcal{L}^+(h)$ fails to be full type
for every single minimal cut system in $V$.  Fix a measure on
$\mathcal{L}^+(h)$.  Fix attention on a single minimal cut system
$C$.  As above let $P$ be the closure of $S \setminus \bdy C$ in the
path metric.  Let $\Gamma = \Gamma(\mathcal{L}^+(h), C)$ be the
Whitehead graph.  By hypothesis there is a disk $D \in C$ giving
$\Gamma$ a cut vertex.

As $D$ gives a cut vertex, there is an essential arc $\gamma$ properly
embedded in $P$ with $\gamma \cap \mathcal{L}^+(h) = \emptyset$.  We
use $\gamma$ to do a {\em disk replacement}: choose a properly
embedded arc $\delta \subset D$ joining the endpoints of $\gamma$.  As
$V \setminus C$ is a ball the curve $\gamma \cup \delta$ bounds a disk $E$.
Let $D', D''$ be the components of $D \setminus \delta$.  One of the
disks $D' \cup E$ or $D'' \cup E$ is nonseparating in $V$.  Thus form
a new minimal cut system $C'$ by removing $D$ and adding this new
disk.

Note that $\iota(C', \mathcal{L}^+(h)) < \iota(C, \mathcal{L}^+(h))$,
by minimality of $\mathcal{L}^+(h)$.  (A cut system is a
measured lamination when given the counting measure.)  Also, by hypothesis,
$\mathcal{L}^+(h)$ again fails to be of full type with respect to
$C'$.  So we may produce a sequence of cut systems, $C^{(n)}$, which
have decreasing intersection number with $\mathcal{L}^+(h)$.  It
follows that $C^{(n)}$ is unbounded in $\ML(S)$.

Choose a sequence $r_n \in \RR_+$ (with $r_n \to 0$) so that the
sequence of measured laminations $\{r_n \cdot C^{(n)}\}$ is bounded in
$\ML(S) \setminus 0$.  Passing to a convergent subsequence let
$\mathcal{L} = \lim (r_n C^{(n)})$.  As $\iota(C^{(n)},
\mathcal{L}^+(h))$ is bounded, $\iota(\mathcal{L}, \mathcal{L}^+(h)) =
0$.  As in the proof of Theorem~\ref{Thm:Hempel}, minimality and
unique ergodicity imply that $\mathcal{L}^+(h) = \mathcal{L}$ and
\Hempel\ does not hold.
\end{proof}

\begin{lemma}
\label{Lem:K'ImpliesK}
If a minimal lamination $\mathcal{L}$ is full type with respect to
some minimal cut system then $\mathcal{L}$ is full type with respect
to some maximal cut system.
\end{lemma}

\begin{proof}
We will prove that if $\mathcal{L}$ is full type with respect to a
non-maximal cut system $C \subset V$ then there is a disk $D \subset
V$ so that $C' = C \cup \{D\}$ is again a cut system and $\mathcal{L}$
remains of full type.

Recall that $P$ is the union of planar surfaces obtained by cutting $S
= \bdy V$ along $C$.  Let $\hat{P}$ be the quotient of $P$ obtained by
identifying each boundary component to a point.  Note that $\Gamma =
\Gamma(\mathcal{L}, C)$ is naturally embedded in $\hat{P}$.  Let $q
\from P \to \hat{P}$ be the quotient map.

Now, a {\em bigon} is any cycle in $\Gamma$ of length two.  A {\em cut
edge} is any edge such that removing the edge, its endpoints, and all
edges adjacent to it, from $\Gamma$ increases the number of components
of $\Gamma$.  Note that if an edge lies on a bigon then it is a cut
edge.  
For every cut edge $E$ which does not lie on a bigon we may add an
extra edge $E' \subset \hat{P}$ to $\Gamma$, with $E' \cap \Gamma =
\bdy E'$, to form a {\em temporary} bigon.  Note also that these extra
edges may be added disjointly, as all simple closed curves in the
sphere separate.

Let $B$ be an innermost bigon (either temporary or not) in $\hat{P}$.
(If there are none set $B = \emptyset$.)  Let $Q$ be the component of
$\hat{P} \setminus B$ which meets no bigons.  Note that none of the
extra edges lie in $Q$.  So choose any edge $F \subset \Gamma$ with
interior contained in the interior of $Q$.  Let $D''$ be a regular
neighborhood of $F$, taken in $\hat{P}$.  Let $D' = q^{-1}(D'')
\subset P$.  Let $\delta = \bdy D' \setminus \bdy P$.  Let $D$ be a
disk in $V$ such that $\bdy D = \delta$.

Let $C' = C \cup \{D\}$ and let $\Gamma' = \Gamma(\mathcal{L}, C')$.
(So $\Gamma'$ is essentially obtained from $\Gamma$ by collapsing the
edge $F$ and removing surplus parallel arcs.)  Now, $\Gamma'$ cannot
have a cut vertex anywhere except at $D$; any such vertex would give
a cut vertex for $\Gamma$.  On the other hand, if $D$ gives a cut vertex
for $\Gamma'$ then the edge $F$ was either part of a bigon (and
$\Gamma'$ has a wave) or $F$ was part of a temporary bigon (and
$\Gamma'$ has a cut vertex without waves.)  This would contradict our
choice of $F$.
\end{proof}

Now \Hempel\ implies \Kobayashi\ by Lemma~\ref{Lem:HImpliesK'} and
Lemma~\ref{Lem:K'ImpliesK}.  This completes the proof of
Theorem~\ref{Thm:KIffH}, and thus of Theorem~\ref{Thm:Main}.


\begin{thebibliography}

\bibitem{Bowditch03}
\textbf{Brian~H Bowditch}, \emph{Tight geodesics in the curve complex},
  \url{http://www.maths.soton.ac.uk/staff/Bowditch/preprints.html}

\bibitem{Bridson99} \textbf{Martin~R Bridson}, \textbf{Andr{\'e}
 Haefliger}, \emph{Metric spaces of non-positive curvature}, Grundlehren
 series 319, Springer--Verlag, Berlin (1999) \MR{1744486}

\bibitem{CassonGordon87}
\textbf{Andrew~J Casson}, \textbf{Cameron~McA Gordon}, \emph{Reducing
  {H}eegaard splittings}, Topology Appl. 27 (1987) 275--283
  \MR{0918537}

\bibitem{CDP90} \textbf{M Coornaert}, \textbf{T Delzant}, \textbf{A
Papadopoulos}, \emph{G\'eom\'etrie et th\'eorie des groupes: Les
groupes hyperboliques de Gromov}, Lecture Notes series 1441,
Springer--Verlag, Berlin (1990) \MR{1075994}

\bibitem{FLP91} \textbf{Fathi}, \textbf{Laudenbach},
\textbf{Po{\'e}naru} (editors), \emph{Travaux de {T}hurston sur les
surfaces}, S\'eminaire Orsay, Ast\'erisque 66--67, Soc. Math.
France, Paris (1979) \MR{0568308}

\bibitem{Gromov87} \textbf{Mikhael Gromov}, \emph{Hyperbolic groups},
from: ``Essays in group theory'', Math. Sci. Res. Inst. Publ 8,
Springer, New York (1987) 75--263 \MR{0919829} 

\bibitem{Hamenstaedt04}
\textbf{U Hamenstaedt}, \emph{Train tracks and the Gromov boundary of the
  complex of curves},
  \arxiv{math.GT/0409611}

\bibitem{Hempel01}
\textbf{John Hempel}, \emph{3--manifolds as viewed from the curve complex},
  Topology 40 (2001) 631--657 \MR{1838999}

\bibitem{Kapovich01}
\textbf{Michael Kapovich}, \emph{Hyperbolic manifolds and discrete groups},
  Progress in Mathematics 183, Birkh\"auser, Boston, MA (2001)
  \MR{1792613}

\bibitem{Klarreich99}
\textbf{Erica Klarreich}, \emph{The boundary at infinity of the curve complex
  and the relative {T}eichm\"{u}ller space},
  \url{http://nasw.org/users/klarreich/research.htm}

\bibitem{Kobayashi88b} \textbf{Tsuyoshi Kobayashi}, \emph{Heights of
simple loops and pseudo-{A}nosov homeomorphisms}, from: ``Braids
(Santa Cruz, CA, 1986)'', Contemp. Math. 78, Amer. Math. Soc.
Providence, RI (1988) 327--338 \MR{0975087}


\bibitem{Masur86}
\textbf{Howard Masur}, \emph{Measured foliations and handlebodies}, Ergodic
  Theory Dynam. Systems 6 (1986) 99--116
  \MR{0837978}

\bibitem{MasurMinsky03}
\textbf{Howard~A Masur}, \textbf{Yair~N Minsky}, \emph{Quasiconvexity in the
  curve complex}
  \arxiv{math.GT/0307083}

\bibitem{MasurMinsky99}
\textbf{Howard~A Masur}, \textbf{Yair~N Minsky}, \emph{Geometry of the complex
  of curves. {I}. {H}yperbolicity}, Invent. Math. 138 (1999) 103--149
  \MR{1714338}

\bibitem{Minsky01} \textbf{Yair~N Minsky}, \emph{Combinatorial and
Geometrical Aspects of Hyperbolic 3--Manifolds}, from: ``Kleinian
groups and hyperbolic 3--manifolds (Warwick, 2001)'', LMS Lecture Note Ser. 299, Cambridge Univ. Press
(2003) 3--40 \MR{2044543}

\bibitem{Shackleton04}
\textbf{Kenneth~J Shackleton}, \emph{{Tightness and computing distances in the
  curve complex}}, \arxiv{math.GT/0412078}

\end{thebibliography}
\end{document}